\documentclass[titlepage,14pt]{extarticle}
\usepackage{amsfonts,amsmath,amssymb,indentfirst}

\usepackage{geometry}
\newcommand{\Li}{\mathop{\rm Li}\nolimits}
\newcommand{\Le}{\mathop{\rm Le}\nolimits}

\def\proof{{\bf Доказательство. }}
\def\hyper#1#2#3#4#5{\vphantom{F}_{#1}F_{#2}\left[{#3\atop#4};#5\right]}

\def\cube#1{[0,1]^{#1}}

\newtheorem{theorem}{Theorem}
\newtheorem{lemma}{Lemma}
\newtheorem{corollary}{Corollary}

\newtheorem{conjecture}{Conjecture}
\newtheorem{claim}{Proposition}
\def\proof{{\bf Proof. }}
\def\qed{\hfill\nobreak$\Box$}
\def\qedInEq{\tag*{$\Box$}}

\def\ze{\widetilde{\zeta}}
\def\inbr#1{ \{ #1 \} }

\def\C{{\mathbb C}}

\def\R{{\mathbb R}}
\def\Q{{\mathbb Q}}
\def\N{{\mathbb N}}

\allowdisplaybreaks[4]

\begin {document}

\begin{center}
{\bf Special Values of Generalized Polylogarithms} \\
{\bf S.A. Zlobin}
\end{center}

\vskip 0.5cm

{\small
We study values of generalized polylogarithms at various points
and relationships among them.
Polylogarithms of small weight at the points 1/2 and -1 are completely
investigated. We formulate
a conjecture about the structure of the linear space generated 
by values of generalized polylogarithms.
}

\vskip 0.5cm

\tableofcontents

\section{Introduction}

Recall that polylogarithms $\Li_k(z)$, $k \in \N$, are defined by
the series $\sum_{n=1}^{\infty} \frac{z^n}{n^k}$. These special functions
are a classical object of study.
At the point $z=1$ they give values of the famous Riemann zeta function, which
is studied in number theory in detail. Recently, interest in
various generalizations of the polylogarithms and the zeta function has 
sharply increased.
In the present paper we consider the functions 
$$
\Li_{s_1, \dots, s_l}(z) = \sum_{n_1 > n_2 > \dots > n_l \ge 1 }
\frac{z^{n_1}}{n_1^{s_1} n_2^{s_2} \dots n_l^{s_l}},
$$
$$
\Le_{s_1, \dots, s_l}(z) = \sum_{n_1 \ge n_2 \ge \dots \ge n_l \ge 1 }
\frac{z^{n_1}}{n_1^{s_1} n_2^{s_2} \dots n_l^{s_l}},
$$
with positive integers $s_j$. The series that define
these generalized polylogarithms have radius of convergence 1.
The generalized polylogarithms themselves can be analytically continued into the domain 
$D = \C \backslash \{ z: |\arg(1-z)| < \pi \}$ (the complex plane
with a cut along the real axis from 1 to $+\infty$).
The series at the point $z=1$ with $s_1>1$ converge and we obtain multiple zeta functions:
$$
\zeta(s_1, s_2, \dots, s_l) = \Li_{s_1, s_2, \dots, s_l} (1),
\quad
\ze(s_1, s_2, \dots, s_l) = \Le_{s_1, s_2, \dots, s_l} (1).
$$
Multiple zeta functions are actively studied
(see for example the survey \cite{ZudilinMZV}). 
Instead, in this paper we consider
generalized polylogarithms at points  $z \ne 1$.
The following values of classical polylogarithms are known:
\begin{equation}
\label{eq1}
\Li_1(z)=-\ln(1-z),
\end{equation}
\begin{equation}
\label{eq2}
\Li_k(-1) = -(1-2^{1-k}) \zeta(k), \quad k > 1,
\end{equation}
\begin{equation}
\label{eq3}
\Li_2 \left( \frac{1}{2} \right)=
\frac{1}{12}\pi^2 - \frac{1}{2} \ln^2 2
=\frac{1}{2} \zeta(2) - \frac{1}{2} \ln^2 2,
\end{equation}
\begin{equation}
\label{eq4}
\Li_3 \left( \frac{1}{2} \right)=
\frac{7}{8} \zeta(3) - \frac{1}{12} \pi^2 \ln 2 
+ \frac{1}{6} \ln^3 2
=\frac{7}{8} \zeta(3) - \frac{1}{2} \zeta(2) \ln 2 
+ \frac{1}{6} \ln^3 2.
\end{equation}
For $k \ge 4$ a similar expression for $\Li_k(1/2)$ 
in terms of $\ln 2$ and multiple zeta values is unknown 
(and presumably does not exist).
Identity (\ref{eq1}) follows from the Taylor series for $\ln(1-z)$.
The equality
$$
\Li_k(1)+\Li_k(-1)=2 \sum_{n=1}^{\infty} \frac{1}{(2n)^k}.
$$
yields (\ref{eq2}) directly.
Formula (\ref{eq3}) was already known to L.~Euler. 
Using the equality
$$
\zeta(2)= 2 \sum_{n=1}^{\infty} \frac{1}{n^2 \cdot 2^n} +\ln^2 2,
$$
he calculated $\zeta(2)=1{.}644934\dots$ 
with six digits precision in the year 1731; the value
$\zeta(2)=\pi^2/6$ was found by Euler in 1735.
The equality (\ref{eq4}) is due to J.~Landen.
Proofs of (\ref{eq3}) and (\ref{eq4}) can be found
in \cite{lewin} (see (1.16) and (6.12)).
In the present paper we consider analogs of (\ref{eq1})--(\ref{eq4})
for generalized polylogarithms, and we also prove
several other equalities.

\section{Integral representations of generalized polylogarithms}

For a vector $\vec{s}=(s_1,\dots,s_l)$ define its {\it weight} 
$w(\vec{s})=s_1+\cdots+s_l$ and its {\it length} $l(\vec{s})=l$.
There exist various integral representations of generalized polylogarithms
$\Li_{\vec{s}}(z)$ and $\Le_{\vec{s}}(z)$
as integrals of dimension equal to the weight of the vector $\vec{s}$.

\begin{theorem}
\label{Int4Polylogarithms}
The following integral representations of generalized polylogarithms 
are valid:
\begin{align*}
\Le_{s_1, s_2, \dots, s_l} (z) & = z \int_{\cube{m}} \frac{dx_1 dx_2
\dots dx_{m}} {\prod_{j=1}^l (1-z x_1 x_2 \dots x_{r_j}) }, \\
\Li_{s_1, s_2, \dots, s_l} (z) & = z^l \int_{\cube{m}}
\frac{\prod_{j=1}^{l-1} (x_1 x_2 \dots x_{r_j}) dx_1 dx_2 \dots dx_{m}}
{\prod_{j=1}^l (1-z x_1 x_2 \dots x_{r_j}) } ,
\end{align*}
where $r_j = s_1 + s_2 + \dots + s_j$ and $m = r_l = w(\vec{s})$.
\end{theorem}
This theorem is proved by expanding each fraction of the form $1/(1-t)$ 
into a geometric series and then integrating termwise
over a cube (a more general proposition was proved in \cite[Lemma 2]{zl5}).
The integrals in Theorem~\ref{Int4Polylogarithms} analytically continue
the generalized polylogarithms into the domain
$D = \C \backslash \{ z: |\arg(1-z)| < \pi \}$.

Generalized polylogarithms 
also have representations in the form of {\it Chen iterated integrals}.
To a vector $\vec{s}=(s_1, s_2, \dots, s_l)$ with positive integer components
we assign the word $x_0^{s_1-1} x_1 \cdots x_0^{s_l-1} x_1$
on the alphabet $\{ x_0, x_1 \}$
(here $x_i^k$ means the letter $x_i$ written $k$ times one after another).
One can recover the vector from a word that ends in $x_1$.
Introduce the following differential forms:
$$
\omega_{x_0}(t)=\frac{dt}{t}, \quad
\omega_{x_1}(t)=\frac{dt}{1-t}, \quad
\widetilde\omega_{x_0}(t)=\frac{dt}{t}, \quad
\widetilde\omega_{x_1}(t)=\frac{dt}{t(1-t)}.
$$
\begin{theorem}
\label{IteratedInt4Polylogarithms}
For any $z \in D$ the following integral representations are valid:
\begin{equation}
\label{LiIteratedInt}
\Li_{\vec{s}}(z)=
\int_{0}^z \omega_{y_1}(t_1) \int_{0}^{t_1} 
\omega_{y_2}(t_2) \cdots \int_{0}^{t_{m-1}} \omega_{y_m}(t_m),
\end{equation}
$$
\Le_{\vec{s}}(z)= 
\int_{0}^z  \widetilde\omega_{y_1}(t_1)
\int_{0}^{t_1}  \widetilde\omega_{y_2}(t_2) \cdots
\int_{0}^{t_{m-2}} \widetilde\omega_{y_{m-1}}(t_{m-1})
\int_{0}^{t_{m-1}} \omega_{y_m}(t_m),
$$
where the word $y_1 y_2 \cdots y_m$ $(y_i \in \{ x_0, x_1 \})$
corresponds to the vector $\vec{s}$ (in particular, $m=w(\vec{s})$).
\end{theorem}
For example, for the vector $\vec{s}=(2,1)$ its corresponding word is
$x_0 x_1 x_1$, and we have
$$
\Li_{2,1}(z)=
\int_{0}^z \frac{1}{t_1} \int_{0}^{t_1} 
\frac{1}{1-t_2} \int_{0}^{t_2} \frac{1}{1-t_3},
$$
$$
\Le_{2,1}(z)=
\int_{0}^z \frac{1}{t_1} \int_{0}^{t_1} 
\frac{1}{t_2(1-t_2)} \int_{0}^{t_2} \frac{1}{1-t_3}.
$$
Theorem \ref{IteratedInt4Polylogarithms} can be obtained
from Theorem \ref{Int4Polylogarithms} by the change of variables
$$x_1=t_1/z, \quad x_2=t_2/t_1, \quad \dots , \quad  x_m=t_m/t_{m-1}.$$
A more direct proof (for $\Li_{\vec{s}}(z)$)
can be found for example in \cite[\S 6]{ZudilinMZV}.

\section{Identities for generalized polylogarithms}

All known relations for values of (generalized) polylogarithms,
in particular (\ref{eq1})-(\ref{eq4}), are weight-homogeneous 
in the following sense. To the series $\Li_{\vec{s}}(z)$
(or $\Le_{\vec{s}}(z)$) with $z \ne 0$ assign the weight equal to the weight
of the vector $\vec{s}$. Thus the series $\zeta(s_1,\dots,s_l)$ also has
weight $w(\vec{s})$, and values of the logarithm have weight
equal to 1, by (\ref{eq1}). Rational numbers by definition have zero weight.
The weight of the product of two quantities is the sum of their weights.

The functions $\Le_{\vec{s}}(z)$ and $\Li_{\vec{s}}(z)$
can each be expressed linearly in terms of the other. Namely
(see \cite{ulansky}):
\begin{equation}
\label{eq5}
\Le_{s_1, \dots, s_l}(z) = \sum_{\vec{p}} \Li_{\vec{p}}(z), \quad
\Li_{s_1, \dots, s_l}(z) = \sum_{\vec{p}} (-1)^{\alpha(\vec{p})}
\Le_{\vec{p}}(z),
\end{equation}
where $\vec{p}$ ranges over all vectors of the form
$(s_1 * s_2 * \dots * s_l)$. Here the symbol
'*' can be either a plus sign '+' or a comma ',', and $\alpha(\vec{p})$
is equal to the number of plus signs '+'. For example,
\begin{align*}
\Le_{s_1,s_2,s_3}(z)&=\Li_{s_1,s_2,s_3}(z) + 
\Li_{s_1 + s_2,s_3}(z) + \Li_{s_1,s_2 + s_3}(z)
+ \Li_{s_1 + s_2 + s_3}(z), \\
\Li_{s_1,s_2,s_3}(z)&=\Le_{s_1,s_2,s_3}(z) -
\Le_{s_1 + s_2,s_3}(z) - \Le_{s_1,s_2 + s_3}(z)
+ \Le_{s_1 + s_2 + s_3}(z).
\end{align*}

The following explicit expression for certain generalized
polylogarithms is known (see, for example, \cite[Section 1]{waldschmidt}):
\begin{equation}
\label{eq6}
\Li_{\inbr{1}_k}(z) = \frac{(- \ln (1-z))^k}{k!}
\end{equation}
($\inbr{a}_k$ is the $k$-dimensional
vector all of whose coordinates are equal to $a$).
Expressions in terms of elementary functions for other generalized polylogarithms are unknown.

\begin{theorem}
If $a$, $b$ are positive integers and $z \in D$, then
$$
\Li_{a,b}(z)= (-1)^b \sum_{k=1}^a \binom{a+b-k-1}{b-1} \Li_{k,a+b-k}(z)
$$
$$
+\sum_{k=1}^b (-1)^{b-k} \binom{a+b-k-1}{a-1} \Li_k(z) \Li_{a+b-k}(z).
$$
\end{theorem}
\proof
It is enough to prove the identity inside the unit disk
$|z|<1$ (then by virtue of analytic continuation it is valid in $D$).
By the definition of the generalized polylogarithm we have
$$
\Li_{a,b}(z)=\sum_{n_1>n_2 \ge 1} \frac{z^{n_1}}{n_1^a n_2^b}=
\sum_{n_1>n_2 \ge 1} \frac{z^{n_1}}{n_1^a (n_1-n_2)^b}.
$$
Expand the last fraction 
into a sum of partial fractions with respect to the variable $n_1$:
$$
\Li_{a,b}(z)=\sum_{n_1>n_2 \ge 1} z^{n_1} \left(
(-1)^b \sum_{k=1}^a \binom{a+b-k-1}{b-1} \frac{1}{n_1^k n_2^{a+b-k}}
\right.
$$
$$
\left.
+\sum_{k=1}^b (-1)^{b-k} \binom{a+b-k-1}{a-1}
\frac{1}{ (n_1-n_2)^{k} n_2^{a+b-k}}
\right).
$$
It is left to note that
$$
\sum_{n_1>n_2 \ge 1} \frac{z^{n_1}}{n_1^k n_2^{a+b-k}}=
\Li_{k,a+b-k}(z)
$$
by definition, and
\begin{equation}
\sum_{n_1>n_2 \ge 1} \frac{z^{n_1}}{ (n_1-n_2)^{k} n_2^{a+b-k}}=
\sum_{m=1}^{\infty} \sum_{n_2=1}^{\infty} 
\frac{z^{m+n_2}}{ m^k n_2^{a+b-k}}=
\Li_k(z) \Li_{a+b-k}(z).
\qedInEq
\end{equation}
\begin{corollary}
\label{cor1}
If $n$ is a positive integer and $z \in D$, then
\begin{equation}
\label{eq7}
2 \Li_{1,2n-1}(z) = 
\sum_{k=1}^{2n-1} (-1)^{k+1}  \Li_k(z) \Li_{2n-k}(z).
\end{equation}
\end{corollary}
\proof
Set $a=1$ and $b=2n-1$ in the theorem. \qed

Equality (\ref{eq7}) was obtained in \cite[Section 6.3]{BradleySpecial}
(at least for $z=1/2$).

\begin{corollary}
\label{cor2}
If $z \in D$, then
$$
2 \Li_{2,1}(z) + \Li_{1,2}(z)= \Li_1(z) \Li_2(z).
$$
\end{corollary}
\proof
Set $a=2$ and $b=1$ in the theorem. \qed
\begin{corollary}
\label{cor3}
If $z \in D$, then
$$
4 \Li_{3,1}(z) + 2 \Li_{2,2}(z)= \Li_2(z)^2.
$$
\end{corollary}
\proof
From the theorem with $a=3$ and $b=1$ it follows that
$$
2 \Li_{3,1}(z)+\Li_{2,2}(z)+\Li_{1,3}(z) = \Li_1(z) \Li_3(z).
$$
Multiplying this equality by 2 and subtracting (\ref{eq7}) with $n=2$,
we obtain the required equality. \qed

Corollaries \ref{cor2} and \ref{cor3} can also be obtained with the help of the
{\it shuffle relations} (see for example \cite[\S 5]{ZudilinMZV}), as follows.
Suppose that $z \in (0,1)$, and that $\vec{s_1}$, $\vec{s_2}$ are vectors with positive
integer components and weights $p=w(\vec{s_1})$, $q=w(\vec{s_2})$.
Formula (\ref{LiIteratedInt}) for $\vec{s_1}$ and $\vec{s_2}$ can be rewritten 
as 
\begin{align*}
\Li_{\vec{s_1}}(z) & = 
\int_{z > t_{1} > t_{2} > \cdots > t_{p} > 0}
\omega_1 \omega_2 \cdots \omega_p, \\
\Li_{\vec{s_2}}(z) & =
\int_{z > t_{p+1} > t_{p+2} > \cdots > t_{p+q} > 0}
\omega_{p+1} \omega_{p+2} \cdots \omega_{p+q},
\end{align*}
where $\omega_i$ is either $dt_i/t_i$ or $dt_i/(1-t_i)$.
Denote by $S'_{p+q}$ the set of permutations $\sigma$ of $\{ 1, \dots, p+q \}$
such that
$$
\sigma^{-1}(1) < \sigma^{-1}(2) < \cdots < \sigma^{-1}(p), \quad
\sigma^{-1}(p+1) < \sigma^{-1}(p+2) < \cdots < \sigma^{-1}(p+q).
$$
Splitting the Cartesian product 
\begin{multline*}
\{ z>t_1>t_2>\cdots>t_p>0 \} \times
\{ z>t_{p+1}>t_{p+2}>\cdots>t_{p+q}>0 \} \\
= \bigsqcup_{\sigma \in S'_{p+q} }
\{ z > t_{\sigma(1)} > t_{\sigma(2)} > \cdots > t_{\sigma(p+q)} > 0 \}
\ \sqcup  \mbox{ \{ a set of measure zero \} },
\end{multline*}
we obtain
$$
\Li_{\vec{s_1}}(z) \Li_{\vec{s_2}}(z) = \sum_{\sigma \in S'_{p+q}}
\int_{z > t_{\sigma(1)} > t_{\sigma(2)} > \cdots > t_{\sigma(p+q)} > 0}
\omega_{\sigma(1)} \omega_{\sigma(2)} \cdots \omega_{\sigma(p+q)}.
$$
Thus the product $\Li_{\vec{s_1}}(z) \Li_{\vec{s_2}}(z)$
is a linear combination of $\Li_{\vec{t}}(z)$ with positive integer
coefficients (the shuffle relation);
the sum of the coefficients is equal to the binomial coefficient
$\binom{w(\vec{s_1})+w(\vec{s_2})}{w(\vec{s_1})}$, and the vectors $\vec{t}$ satisfy the equalities
$$
l(\vec{t})=l(\vec{s_1})+l(\vec{s_2}), \quad
w(\vec{t})=w(\vec{s_1})+w(\vec{s_2}).
$$
Moreover this identity  holds in the domain $D$ by analytic continuation.
Thus the second proof of Corollaries \ref{cor2} and \ref{cor3} is complete.

In what follows we will need the following shuffle relations:
\begin{align}
\label{ShuffleWeight5Eq1}
\Li_1(z) \Li_4(z) &=
2 \Li_{4,1}(z)+\Li_{3,2}(z)+\Li_{2,3}(z)+\Li_{1,4}(z), \\
\label{ShuffleWeight5Eq2}
\Li_2(z) \Li_3(z) &=
6 \Li_{4,1}(z)+3 \Li_{3,2}(z)+\Li_{2,3}(z), \\
\label{ShuffleWeight5Eq3}
\Li_1(z) \Li_{2,2}(z) &=
2 \Li_{2,2,1}(z)+ 2 \Li_{2,1,2}(z)+\Li_{1,2,2}(z).
\end{align}

Suppose we have a vector $\vec{s}=(s_1, s_2, \dots, s_l)$ with positive
integer components. As before, to the vector
$\vec{s}$ assign the word $x_0^{s_1-1}x_1\cdots x_0^{s_l-1}x_1=v x_1$.
Let \(\sigma\ \) be the mapping acting on these words which interchanges
the letters  \( x_0 \) and \( x_1 \).
Define the dual vector $\vec{s'}$ be the vector that corresponds
to the word $\sigma (v)x_1$. For example, if $\vec{s}=(2,1,3)$, then its
corresponding word is $x_0 x_1^2 x_0^2 x_1 = v x_1$, and
$\sigma(v) x_1 = x_1 x_0^2 x_1^3$, $\vec{s'}=(1,3,1,1)$.
The following theorem is valid (see \cite[Lemma 12]{zl5}).

\begin{theorem}
\label{dualityTh}
(Duality)
If $z \in D$, then
$$
\Le_{\vec{s}} \left( \frac{-z}{1-z} \right) = - \Le_{\vec{s'}} (z).
$$
\end{theorem}
\begin{corollary}
\label{FullSummCor}
(On the sum of polylogarithms of constant weight)
Let $W_n$ be the set of all vectors having weight $n$.
Then for any $z \in D$ the following equalities hold:
$$
\Le_{\inbr{1}_n} (z)=
\sum_{\vec{p} \in W_n} \Li_{\vec{p}}(z)=-\Li_{n} \left( \frac{-z}{1-z} \right),
$$
$$
\Li_{\inbr{1}_n} (z)=
\sum_{\vec{p} \in W_n} (-1)^{n-l(\vec{p})} \Le_{\vec{p}}(z)
= \frac{(- \ln (1-z))^k}{k!}.
$$
\end{corollary}
\proof
Formula (\ref{eq5}) yields identities
$$
\sum_{\vec{p} \in W_n} \Li_{\vec{p}}(z)= \Le_{\inbr{1}_n}(z), \quad
\sum_{\vec{p} \in W_n} (-1)^{n-l(\vec{p})} \Le_{\vec{p}}(z)
= \Li_{\inbr{1}_n}(z).
$$
Theorem \ref{dualityTh} with $\vec{s}=(n)$ leads to
$$
\Le_{\inbr{1}_n}(z)=-\Li_{n} \left( \frac{-z}{1-z} \right),
$$
which proves the first required pair of equalities.
The second pair follows from identity (\ref{eq6}). \qed

The generalized polylogarithm $\Li_{\vec{s}}(z)$ is
a linear combination of $\Le_{\vec{t}}(z)$, so by Theorem \ref{dualityTh}
the function $\Li_{\vec{s}}(\frac{-z}{1-z})$
is a linear combination of $\Le_{\vec{t'}}(z)$, 
and hence of $\Li_{\vec{u}}(z)$. 
The explicit formula (see \cite[Theorem 1]{ulansky})
is the following. (In the statement, instead of vectors for polylogarithms, the
corresponding words are used.)
\begin{theorem}
\label{LiTransformation}
If $z \in D$ and $s_1$, \dots, $s_l$ are positive integers, then 
$$
\Li_{x_0^{s_1-1} x_1 \cdots x_0^{s_l-1} x_1} \left( \frac{-z}{1-z} \right) =
(-1)^l \sum_{w(p_1)=s_1-1, \dots, w(p_l)=s_l-1}
\Li_{p_1 x_1 \cdots p_l x_1} (z).
$$
\end{theorem}
For example, if $s_1=3$, $s_2=1$, then 
\begin{align*}
\Li_{3, 1} \left( \frac{-z}{1-z} \right) & =
\sum_{w(p)=2} \Li_{p x_1 x_1} (z) \\
& = \Li_{x_0 x_0 x_1 x_1} (z) + \Li_{x_0 x_1 x_1 x_1} (z) 
+ \Li_{x_1 x_0 x_1 x_1} (z) + \Li_{x_1 x_1 x_1 x_1}(z) \\
& = \Li_{3, 1} (z) + \Li_{2, 1, 1} (z) 
+ \Li_{1, 2, 1} (z) + \Li_{1, 1, 1, 1} (z).
\end{align*}
  
In \cite[Section 7]{BradleySpecial} an identity was obtained linking
generalized polylogarithms at the points $z$ and $1-z$ (in that paper
it was called a {\it H\"older convolution}, since the authors worked
with an argument inversely proportional to ours).
\begin{theorem}
\label{HolderTh}
Suppose $m \ge 2$ and $y_1 y_2 \cdots y_m$ is a word on the alphabet
$\{ x_0, x_1 \}$, where $y_1=x_0$ and $y_m=x_1$. Denote by $\tau$ the function
explicitly defined on this alphabet by $\tau(x_0)=x_1$, $\tau(x_1)=x_0$.
Then for any $z \in \C$ not lying on either of the rays
$\{ z \in \R: z \le 0  \}$, $\{ z \in \R: z \ge 1 \}$, we have
$$
\zeta(y_1 y_2 \cdots y_m)=\sum_{i=1}^{m+1} 
\Li_{y_i \cdots y_m}(z) \Li_{\tau(y_{i-1}) \cdots \tau(y_1)}(1-z) ,
$$
where we suppose $\Li_{\emptyset}(z) \equiv 1$.
\end{theorem}
For example, applying Theorem \ref{HolderTh} to the word $x_0^{m-1} x_1$, 
we obtain
\begin{align*}
\zeta(m) &=\sum_{i=1}^m \Li_{x_0^{m-i} x_1}(z) \Li_{x_1^{i-1}}(1-z)
+ \Li_{x_0 x_1^{m-1}}(1-z) \\
&=\sum_{i=1}^m \Li_{m-i+1}(z) \Li_{\inbr{1}_{i-1}}(1-z)
+ \Li_{2, \inbr{1}_{m-2}}(1-z).
\end{align*}
This theorem was proved in \cite{BradleySpecial} for $z \in (0,1)$ with the
help of the representation (\ref{LiIteratedInt}). In the broader domain
the identity holds by analytic continuation.

In \cite[Theorem 2]{ulansky} an explicit formula is given which 
expresses $\Li_{\vec{s}}(1-z)$ as a polynomial with rational coefficients
in $\Li_{\vec{t}}(z)$, $\ln z$, and multiple zeta values.
Another approach to getting formulas for a generalized polylogarithm with linear-fractional transformation of its argument is developed in \cite{MinhGenerating}. Using a formal generating series for polylogarithms, the transformations
$z \to 1-z$, $z \to 1-1/z$, and $z \to 1/z$ are studied in this work.

\section{The alternating multiple zeta function}
Define the {\it alternating multiple zeta function} 
$$
\zeta(s_1,\dots,s_l; \sigma_1, \dots, \sigma_l)
=
\sum_{n_1 > \cdots > n_l \ge 1} \frac{\sigma_1^{n_1} \cdots \sigma_l^{n_l}}
{n_1^{s_1} \cdots n_l^{s_l}}
$$
for positive integers $s_j$ and $\sigma_j=\pm 1$
(if all $\sigma_j$ are complex roots of unity, the name
{\it colored multiple zeta function} is also used).
Sometimes these series are called (alternating) {\it Euler sums}.
Often the exponent $s_j$ and the sign $\sigma_j$ are combined into one symbol:
$s_j$ if $\sigma_j=1$ and $\overline s_j$ if $\sigma_j=-1$.
For example,
$\zeta(\overline s_1, s_2, \dots, s_l) = \Li_{s_1, s_2, \dots, s_l}(-1)$.
The alternating multiple zeta function 
is surely deserving of individual study, but in the 
present paper it plays an auxiliary role. 

We will also need the function
$$
\mu(a_1, \dots,  a_l) = 
\sum_{k_1=1}^{\infty} \cdots \sum_{k_l=1}^{\infty}
\frac{a_1^{k_1} \cdots a_l^{k_l}}
{(k_1+k_2+\dots+k_l) \cdots (k_{l-1}+k_l) k_l}.
$$
If $a_1=-1$ and $a_j=\pm 1$ when $j>1$, then the series converges and
with the help of the change of variables $n_j=k_j+\dots+k_l$, we obtain
\begin{equation}
\label{EulerSumToAlterZeta}
\mu(a_1, \dots,  a_l)=\zeta(\inbr{1}_l; \sigma_1, \dots, \sigma_l),
\end{equation}
where $\sigma_1=-1$ and $\sigma_j \equiv a_{j-1}+a_j-1 \pmod 4$ for $j>1$.

\begin{theorem}
\label{ToEulerSumsTh}
If $s_j$ are positive integers and $z \in \C$, $|z|<1$, then
$$
\Li_{s_1, \dots, s_l} \left( \frac{-z}{1-z} \right)=
(-1)^l  \mu(z, \inbr{1}_{s_l-1}, \dots, z, \inbr{1}_{s_1-1}).
$$
\end{theorem}
\proof
It is enough to prove the theorem for $z \in (-1,0)$ (by virtue of analytic continuation).
Denote by $y_1 y_2 \cdots y_m$ the word
$x_0^{s_1-1} x_1 \cdots x_0^{s_l-1} x_1$. 
By Theorem \ref{IteratedInt4Polylogarithms} we have
$$
\Li_{s_1, \dots, s_l} \left( \frac{-z}{1-z} \right)=
\int_{0}^{\frac{-z}{1-z}} \omega_{y_1}(t_1) \int_{0}^{t_1} 
\omega_{y_2}(t_2) \cdots \int_{0}^{t_{m-1}} \omega_{y_m}(t_m).
$$
After the change of variables $t_{m+1-i}=-z/(1-z) \cdot (1-u_i)$ we obtain
$$
\Li_{s_1, \dots, s_l} \left( \frac{-z}{1-z} \right)=
\int_{0}^{1} \nu_{y_m}(u_1) \int_{0}^{u_1} 
\nu_{y_{m-1}}(u_2) \cdots \int_{0}^{u_{m-1}} \nu_{y_1}(u_m),
$$
where
$$
\nu_{x_0}(u)=\frac{du}{1-u}, \quad \nu_{x_1}(u)=\frac{-z \, du}{1-zu}.
$$
Successively integrating over the variables $u_m$, $u_{m-1}$, \dots, $u_1$,
we arrive at the expression
\begin{equation*}
(-1)^l  \mu(z, \inbr{1}_{s_l-1}, \dots, z, \inbr{1}_{s_1-1}). \qedInEq
\end{equation*}
\begin{corollary}
\label{ToEulerSumsCor}
For any positive integers $s_j$ the following equality holds
$$
\Li_{s_1, \dots, s_l} \left( \frac{1}{2} \right)=
(-1)^l  \mu(-1, \inbr{1}_{s_l-1}, \dots, -1, \inbr{1}_{s_1-1}).
$$
\end{corollary}
\proof
In Theorem \ref{ToEulerSumsTh} take the limit as $z \to -1$. \qed
\begin{corollary}
\label{MuByLiCor}
If $s_j$ are positive integers and $|z|<1$, then
$$
\mu(z, \inbr{1}_{s_l-1}, \dots, z, \inbr{1}_{s_1-1})=
\sum_{w(p_1)=s_1-1, \dots, w(p_l)=s_l-1}
\Li_{p_1 x_1 \cdots p_l x_1} (z).
$$
\end{corollary}
\proof
The required identity directly follows from
Theorems \ref{LiTransformation} and \ref{ToEulerSumsTh}. \qed

Corollary \ref{ToEulerSumsCor} was proved in \cite[(6.8)]{BradleySpecial}.
It determines a one-to-one correspondence between
values of polylogarithms $\Li_{\vec{s}}(1/2)$ and values of
the alternating multiple zeta function with $s_j=1$, $\sigma_1=-1$
and arbitrary $\sigma_j=\pm 1$ at $j>1$.
Corollary \ref{MuByLiCor} for $z=-1$ (which can be obtained 
by a limiting process) was proved in \cite[Theorem 9.4]{BradleySpecial}.

On the basis of experimental data we advance the following conjecture. 

\begin{conjecture}
\label{EulerSumsConj}
Any value of the alternating multiple zeta function of weight $w$
is a rational linear
combination of numbers of any one of the following forms: \\
1) $\Li_{\vec{s}}(1/2)$, $w(\vec{s})=w$ \\
2) $\Le_{\vec{s}}(1/2)$, $w(\vec{s})=w$ \\
3) $\Li_{\vec{s}}(-1)$, $w(\vec{s})=w$ \\
4) $\Le_{\vec{s}}(-1)$, $w(\vec{s})=w$ \\
5) or $\zeta(\inbr{1}_w; -1, \sigma_2, \dots, \sigma_w)$
with arbitrary $\sigma_j=\pm 1$ at $j>1$.
\end{conjecture}
We show that items 1-5 are equivalent. Indeed 
$1 \Leftrightarrow 2$ and $3 \Leftrightarrow 4$ by (\ref{eq5}),
$2 \Leftrightarrow 4$ by Theorem \ref{dualityTh},
$1 \Leftrightarrow 5$ by Corollary \ref{ToEulerSumsCor}.
In \cite[Section 7]{BradleySpecial} it was shown that the statement of 
Conjecture \ref{EulerSumsConj} is true for the (nonalternating)
multiple zeta function. For that it is enough to make use of Theorem
\ref{HolderTh} for $z=1/2$, and then to apply the shuffle relation 
for every product of two polylogarithms. Also it is not difficult
to prove the statement for the alternating multiple zeta function of length 
$\le 2$. Additionally the conjecture was checked by the author for weight
$\le 8$. The verification was implemented by computing
Euler sums with high precision (see \cite[Section 7]{BradleySpecial})
and using the PSLQ algorithm (see \cite{PSLQ}) for finding linear relations
with integer coefficients.

\section{Values of generalized polylogarithms at the point $z=1/2$}
\label{HalfSection}

\begin{theorem}
\label{th2}
(On the sum of polylogarithms of constant weight at the point $z=1/2$)
Let $W_n$ be the set of all vectors having weight $n$.
Then the following equalities hold:
$$
\Le_{\inbr{1}_n} \left( \frac{1}{2} \right) =
\sum_{\vec{p} \in W_n} \Li_{\vec{p}} \left( \frac{1}{2} \right)
= (1-2^{1-k}) \zeta(k),
$$
$$
\Li_{\inbr{1}_n} \left( \frac{1}{2} \right) =
\sum_{\vec{p} \in W_n} (-1)^{n-l(\vec{p})} \Le_{\vec{p}}
\left( \frac{1}{2} \right)
= \frac{\ln^n 2}{n!}.
$$
\end{theorem}
\proof
The required equalities follow from Corollary \ref{FullSummCor}
with $z=1/2$ together with formula (\ref{eq2}). \qed

Another relation for values of generalized polylogarithms at the point $z=1/2$
comes from Theorem \ref{HolderTh} applied to $z=1/2$.

From (\ref{eq7}) we obtain the value $\Li_{1,2n-1}(1/2)$ for any positive
integer $n$; this value lies in the ring
$\Q[\Li_1(1/2),\Li_2(1/2),\dots,\Li_{2n}(1/2)]$.

\begin{lemma}
\label{lemma1}
The following integral representation is valid:
$$
\Li_{\inbr{1}_m,2,\inbr{1}_n} \left( \frac{1}{2} \right)=
\frac{1}{m! (n+1)!} \int_{1/2}^1 (\ln(2t))^m (-\ln t)^{n+1} \frac{dt}{1-t}.
$$
\end{lemma}
\proof
The integral representation (\ref{LiIteratedInt}) yields
\begin{align*}
\Li_{\inbr{1}_m,2,\inbr{1}_n} \left( \frac{1}{2} \right)= &
\int_0^{1/2} 
\left( \int_{t_{m+1}}^{1/2} \frac{dt_{m}}{1-t_{m}} \cdots 
\int_{t_{2}}^{1/2} \frac{dt_{1}}{1-t_{1}} \right) \\
& \times
\left( \int_0^{t_{m+1}} \frac{dt_{m+2}}{1-t_{m+2}} \cdots 
\int_0^{t_{m+n+2}} \frac{dt_{m+n+2}}{1-t_{m+n+2}} \right)
\times
\frac{dt_{m+1}}{t_{m+1}} \\
= & \int_0^{1/2} \frac{(\ln(2(1-t_{m+1}))^{m}}{m!} \cdot
\frac{(-\ln(1-t_{m+1}))^{n+1}}{(n+1)!} \cdot
\frac{dt_{m+1}}{t_{m+1}}.
\end{align*}
The substitution $t=1-t_{m+1}$ completes the proof of the lemma. \qed
\begin{corollary}
For any integers $m \ge 0$ and $n \ge 0$ the following
equalities hold: 
\begin{equation}
\label{eq8}
\Li_{\inbr{1}_m,2,\inbr{1}_n} \left( \frac{1}{2} \right)=
\sum_{k=0}^m (-1)^k \binom{n+1+k}{n+1}
\Li_{2,\inbr{1}_{n+k}} \left( \frac{1}{2} \right)
\frac{\ln^{m-k} 2}{(m-k)!},
\end{equation}
\begin{equation*}
\Li_{\inbr{1}_m,2,\inbr{1}_n} \left( \frac{1}{2} \right)=
\frac{1}{n+1} \sum_{k=0}^n (-1)^k  \binom{m+k}{m}
\Li_{\inbr{1}_{m+k},2} \left( \frac{1}{2} \right)
\frac{\ln^{n-k} 2}{(n-k)!}.
\end{equation*}
\end{corollary}
\proof
The first equality follows by the substituting the identity
$$
(\ln(2t))^m=\sum_{k=0}^m (-1)^k \binom{m}{k} (-\ln t)^k \ln^{m-k} 2,
$$
in the integrand of Lemma \ref{lemma1}, and the second by the identity
\begin{equation*}
(-\ln t)^{n}=\sum_{k=0}^{n} (-1)^k \binom{n}{k}
(\ln (2t))^k \ln^{n-k} 2. \qedInEq
\end{equation*}

Lemma \ref{lemma1} was used in \cite{BradleySpecial} to calculate the value 
$\Li_{\inbr{1}_m,2,\inbr{1}_n}$, as
in the following theorem; we propose a new proof of this 
theorem. 

\begin{theorem}
\label{Li121_Th}
For any integers $m \ge 0$ and $n \ge 0$, the value
$\Li_{\inbr{1}_m,2,\inbr{1}_n} (1/2)$ can be expressed
in the form of a polynomial of several variables with rational
coefficients in $\ln2$, $\zeta(a)$, and $\Li_{b}(1/2)$, with
$a \in [n+2, n+m+2]$ and $b \in [m+1, m+n+2]$.
Namely, the following formula is valid:
$$
\Li_{\inbr{1}_m,2,\inbr{1}_n} \left( \frac{1}{2} \right) =
\sum_{k=0}^m (-1)^k \binom{n+1+k}{n+1} \zeta(n+2+k) 
\frac{\ln^{m-k} 2}{(m-k)!}
$$
$$
+
(-1)^{m+1} \sum_{k=0}^{n+1} \binom{m+k}{m}
\Li_{m+1+k} \left( \frac{1}{2} \right)
\frac{\ln^{n+1-k} 2}{(n+1-k)!}.
$$
\end{theorem}
\proof
We first prove the theorem for $m=0$, i.e., the equality
\begin{equation}
\label{Li_1m_2_1n_case_m0}
\Li_{2,\inbr{1}_n} \left( \frac{1}{2} \right) =
\zeta(n+2) 
- \sum_{k=0}^{n+1} \Li_{k+1} \left( \frac{1}{2} \right)
\frac{\ln^{n+1-k} 2}{(n+1-k)!}.
\end{equation}
Using Lemma \ref{lemma1}, we obtain
$$
\Li_{2,\inbr{1}_n} \left( \frac{1}{2} \right) =
\frac{1}{(n+1)!} \int_{1/2}^1 (-\ln t)^{n+1} \frac{dt}{1-t}.
$$
A similar integral was evaluated in \cite{SondowZlobin} (see the proof of 
Theorem 7); we use the same method.
Expanding $1/(1-t)$ into a geometric series gives
\begin{equation}
\label{Li2_1n_eq}
\Li_{2,\inbr{1}_n} \left( \frac{1}{2} \right) =
\frac{1}{(n+1)!} \sum_{r=1}^{\infty} \int_{1/2}^1 t^{r-1} (-\ln t)^{n+1} dt.
\end{equation}
Now consider the equality
$$
\int_{1/2}^1 t^{\sigma+r-1} dt = \frac{1}{\sigma+r}-
\frac{1}{2^{\sigma+r}(\sigma+r)}.
$$
If we differentiate $(n+1)$ times with respect to $\sigma$,
and then set $\sigma=0$, we get
$$
\int_{1/2}^1 t^{r-1} (-\ln t)^{n+1} dt=
(n+1)! \left( \frac{1}{r^{n+2}} - \frac{1}{2^r}
\sum_{k=0}^{n+1} \frac{\ln^{n+1-k} 2}{r^{k+1} (n+1-k)!} \right).
$$
Putting this value of the integral in $(\ref{Li2_1n_eq})$, we obtain
(\ref{Li_1m_2_1n_case_m0}).
(One can also get (\ref{Li_1m_2_1n_case_m0}) from Theorem \ref{HolderTh} by
setting $z=1/2$, $m=n+2$, $y_1 \cdots y_m = x_0^{n+1} x_1$.)

We now prove the theorem for an arbitrary $m$.
In (\ref{eq8}) to each  $\Li_{2,\inbr{1}_{n+k}} \left( \frac{1}{2} \right)$
apply (\ref{Li_1m_2_1n_case_m0}):
\begin{align*}
\Li_{\inbr{1}_m,2,\inbr{1}_n} \left( \frac{1}{2} \right)= &
\sum_{k=0}^m (-1)^k \binom{n+1+k}{n+1}
\frac{\ln^{m-k} 2}{(m-k)!} \\
& \times
\left( \zeta(n+2+k) 
- \sum_{l=0}^{n+1+k} \Li_{l+1} \left( \frac{1}{2} \right)
\frac{\ln^{n+1+k-l} 2}{(n+1+k-l)!}
\right) \\
= & \sum_{k=0}^m (-1)^k \binom{n+1+k}{n+1}
\zeta(n+2+k) \frac{\ln^{m-k} 2}{(m-k)!}  \\
& -\sum_{l=0}^{n+m+1} S_{m,n,l} \Li_{l+1} \left( \frac{1}{2} \right)
\ln^{n+m+1-l} 2 ,
\end{align*}
where 
$$
S_{m,n,l}=\sum_{k=\max(0, l-n-1)}^m 
(-1)^k \binom{n+1+k}{n+1} \frac{1}{(m-k)! (n+1+k-l)!}.
$$
The theorem will be proved if we show that
$S_{m,n,l}=0$ when $l<m$ 
and that
$$
S_{m,n,l}=(-1)^{m} \binom{l}{m} \frac{1}{(n+m+1-l)!}
$$
when $m \le l \le m+n+1$.
If $l \le n+1$ rewrite $S_{m,n,l}$ in the form of 
the hypergeometric function
$$
\frac{1}{m! (n+1-l)!} \cdot \hyper{2}{1}{-m,n+2}{n+2-l}{1},
$$
and if $l>n+1$ in the form
$$
(-1)^{l-n-1} \binom{l}{n+1}
\frac{1}{(m-l+n+1)!} \cdot \hyper{2}{1}{-m+l-n-1,l+1}{l-n}{1}.
$$
In both cases the value of the hypergeometric function
can be evaluated by Gauss' theorem (see \cite[(2.1.3(14)]{Bateman}),
which in this situation is Vandermonde's theorem
$$
\hyper{2}{1}{-k,b}{c}{1}=\frac{(c-b)_k}{(c)_k}, \quad
k=0,1,2,\dots,
$$
where $(a)_n=\Gamma(a+n)/\Gamma(a)$ is the Pochhammer symbol.
We obtain the expressions
$$
\frac{1}{m! (n+1-l)!} \cdot \frac{(-l)_m}{(n+2-l)_m},
$$
$$
(-1)^{l-n-1} \binom{l}{n+1}
\frac{1}{(m-l+n+1)!} \cdot
\frac{(-n-1)_{m-l+n+1}}{(l-n)_{m-l+n+1}},
$$
which can be transformed into the required form. \qed

In Section \ref{Minus1Section} we will show that the values
$\Le_{2, \inbr{1}_n} (1/2)$ for any $n$, and \linebreak
$\Le_{\inbr{1}_m, 2, \inbr{1}_n} (1/2)$ for odd $m+n$,
lie in the ring $\Q[\ln2, \zeta(2), \zeta(3), \dots, ]$.

Theorem \ref{Li121_Th} was proved in \cite[Theorem 8.5]{BradleySpecial}
with the help of generating functions. The method of 
generating functions often turns out to be useful for 
computing values of polylogarithms whose vectors 
are given by a periodic rule.
\begin{theorem}
\label{Li2_Th}
For any nonnegative integer $m$ the values
$\Li_{\inbr{2}_m} (1/2)$ and \linebreak
$\Li_{1,\inbr{2}_m} (1/2)$
belong to the ring $\Q[\ln2, \zeta(2), \zeta(3), \dots]$.
The following identity is valid for any $z \in \C$:
\begin{multline*}
G(z):=\sum_{m=0}^{\infty}
(-1)^m \Li_{\inbr{2}_m} \left( \frac{1}{2} \right) z^{2m} 
+ \sum_{m=0}^{\infty}
(-1)^{m+1} \Li_{1,\inbr{2}_m} \left( \frac{1}{2} \right) z^{2m+1} \\
=
\frac{\Gamma(1/2)}{\Gamma(1+z/2) \Gamma(1/2-z/2)}.
\end{multline*}
\end{theorem}
\proof
We first estimate $\Li_{\inbr{2}_m} (1/2)$ and $\Li_{1,\inbr{2}_m} (1/2)$.
It is clear that
$$
\Li_{\inbr{2}_{m+1}} \left( \frac{1}{2} \right) <
\Li_{1,\inbr{2}_m} \left( \frac{1}{2} \right) =
\sum_{n_1=m+1}^{\infty} \frac{1}{2^{n_1} n_1}
\sum_{n_2=m}^{n_1-1} \frac{1}{n_2^2} \cdots
\sum_{n_m=2}^{n_{m-1}-1} \frac{1}{n_m^2}
\sum_{n_{m+1}=1}^{n_{m}-1} \frac{1}{n_{m+1}^2}.
$$
Using the inequalities
$$
\sum_{n=1}^{r} \frac{1}{n^2} < \sum_{n=1}^{\infty} \frac{1}{n^2} = \zeta(2),
\quad
\sum_{n=p}^{r} \frac{1}{n^2} < \int_{p-1}^{\infty} \frac{dt}{t^2}
=\frac{1}{p-1}, \; p > 1,
$$
we obtain
$$
\Li_{\inbr{2}_{m+1}} \left( \frac{1}{2} \right) <
\Li_{1,\inbr{2}_m} \left( \frac{1}{2} \right) < \frac{\zeta(2)}{(m-1)!}
\sum_{n_1=m+1}^{\infty} \frac{1}{2^{n_1} n_1} <
\frac{\zeta(2) \ln 2}{(m-1)!}.
$$                                                              
Hence both power series have infinite radius of convergence.

By Corollary \ref{ToEulerSumsCor} and relation (\ref{EulerSumToAlterZeta})
we have
$$
\Li_{\inbr{2}_m} \left( \frac{1}{2} \right)=
(-1)^m \mu(\inbr{-1,1}_m)= (-1)^m \zeta( \inbr{\overline 1}_{2m}),
$$
$$
\Li_{1,\inbr{2}_m} \left( \frac{1}{2} \right)=
(-1)^{m+1} \mu(\inbr{-1,1}_m,-1) 
= (-1)^{m+1} \zeta( \inbr{\overline 1}_{2m+1}).
$$
Therefore
$$
G(z) = \sum_{m=0}^{\infty} z^m \sum_{n_1 > n_2 > \cdots > n_{m} \ge 1} 
\frac{ (-1)^{n_1+n_2+\dots+n_{m}} } {n_1 n_2 \cdots n_{m}}
=\prod_{j=1}^{\infty} \left( 1 + (-1)^j \frac{z}{j} \right).
$$
Hence we obtain (see \cite[1.3(1)]{Bateman}) that
$$
G(z)=\frac{\Gamma(1/2)}{\Gamma(1+z/2) \Gamma(1/2-z/2)}.
$$
Then for $|z|<1$ we have $G(z) \ne 0$ and
\begin{align*}
(\ln G(z))' & = \sum_{j=1}^{\infty} (\ln ( 1 + (-1)^j z/j  ) )'=
\sum_{j=1}^{\infty} \frac{1}{1 + (-1)^j z/j } \frac{(-1)^j}{j} \\
& =\sum_{j=1}^{\infty} \sum_{k=0}^{\infty} \frac{(-1)^{jk+j+k} z^{k}}{j^{k+1}}
=\sum_{k=0}^{\infty} (-1)^k \Li_{k+1}((-1)^{k+1}) z^k,
\end{align*}
or 
$$
G(z) = C \exp \left( \sum_{k=1}^{\infty}
\frac{(-1)^{k-1} \Li_{k}((-1)^{k}) z^k}{k} \right).
$$
Moreover $C=1$ from $G(0)=1$. The formula for $G(z)$ allows us to evaluate the coefficient
of $z^k$ in its Taylor series at the point $z=0$.
The equality
\begin{equation*}
\Li_{k}((-1)^{k}) =
\begin{cases}
-\ln 2, & \text{$k$=1}, \\
\zeta(k), & \text{$k$ is even}, \\
-(1-2^{1-k}) \zeta(k), & \text{$k$ is odd $> 1$}, \\
\end{cases}
\end{equation*}
implies that the coefficient lies in the ring
$\Q[\ln2, \zeta(2), \zeta(3), \dots, \zeta(k)]$.
Now the theorem is completely proved. \qed

Theorem \ref{Li2_Th} (for the alternating multiple zeta function)
was implicitly proved in \cite{BradleyEval} (see formulas (13) and (62)).
In \cite{BradleyResolutions} it is explicitly proved
(see Theorems 3 and 4, and also the remarks at the end of Sections
4.2 and 4.4) that the values
$$
\Li_{\inbr{3,1}_n} \left( \frac{1}{2} \right), \quad
\Li_{2,1,\inbr{3,1}_n} \left( \frac{1}{2} \right), \quad
\Li_{1,1,\inbr{3,1}_n} \left( \frac{1}{2} \right), \quad
\Li_{1,\inbr{3,1}_n} \left( \frac{1}{2} \right),
$$
$$
\Li_{\inbr{1,3}_n} \left( \frac{1}{2} \right), \quad
\Li_{3,\inbr{1,3}_n} \left( \frac{1}{2} \right), \quad
\Li_{2,\inbr{1,3}_n} \left( \frac{1}{2} \right), \quad
\Li_{1,\inbr{1,3}_n} \left( \frac{1}{2} \right)
$$
lie in the ring $\Q[\ln2, \zeta(2), \zeta(3), \dots]$.
Moreover, the generating functions are given, from which one can calculate
all those values.

\section{Values of generalized polylogarithms at the point $z=-1$}
\label{Minus1Section}

Values of generalized polylogarithms at the point $z=-1$ can be 
reduced to values at the point $z=1/2$ with the help of Theorems 
\ref{dualityTh} and \ref{LiTransformation}. Thus a conditionally convergent
series of polylogarithms (at $z=-1$) can be equated to the sum of
a more rapidly convergent series. Also there are several specific theorems like those for values at the point $z=1/2$.

\begin{theorem}
(On the sum of polylogarithms of constant weight at the point $z=-1$)
Let $W_n$ be the set of all vectors having weight $n$.
Then the following equalities hold:
$$
\Le_{\inbr{1}_n} (-1) =
\sum_{\vec{p} \in W_n} \Li_{\vec{p}}(-1)=
-\Li_{n} \left( \frac{1}{2} \right)
$$
$$
\Li_{\inbr{1}_n} (-1) =
\sum_{\vec{p} \in W_n} (-1)^{n-l(\vec{p})} \Le_{\vec{p}}(-1)
= \frac{(-\ln 2)^n}{n!}.
$$
\end{theorem}
\proof
The required equalities follow from Corollary \ref{FullSummCor}
with $z=-1$. \qed

Theorem \ref{LiTransformation} yields
$$
\Li_{\inbr{1}_m,2,\inbr{1}_n}(-1)= (-1)^{m+n+1} \left(
\Li_{\inbr{1}_m,2,\inbr{1}_n} \left( \frac{1}{2} \right) +
\Li_{\inbr{1}_{m+n+2}} \left( \frac{1}{2} \right)
\right).
$$
Applying (\ref{eq6}) and Theorem \ref{Li121_Th}, we obtain an explicit formula for
$\Li_{\inbr{1}_m,2,\inbr{1}_n}(-1)$.

In \cite[Theorem 7.2]{flajolet} the following theorem was proved by evaluating complex integrals.
\begin{theorem}
\label{Le_mn_Th}
For any positive integers  $m$, $n$ with odd sum $m+n$
the value $\Le_{m, n}(-1)$ belongs to the ring
$\Q[\ln2, \zeta(2), \dots, \zeta(m+n)]$,
namely, the formula
\begin{multline*}
2 \Le_{m, n}(-1) =
(1-(-1)^n) \zeta(n) \zeta(\overline m)
+\zeta( \overline {m+n} ) \\
+2 (-1)^n \sum_{j+2k=n} \binom{m+j-1}{m-1}  \zeta(\overline {m+j})
	\zeta( \overline {2k} ) \\
+2 (-1)^n \sum_{i+2k=m} \binom{n+i-1}{n-1} \zeta(n+i) \zeta(\overline {2k}),
\end{multline*}
is valid, where $i$, $j$, $k$ are nonnegative integers and we define 
$\zeta(\overline 0)=-1/2$.
If $n=1$, then the following rule applies:
on the right-hand side of the equality the first summand
$2 \zeta(1) \zeta(\overline m)$ cancels
the term in the sum over $i$ with $i=0$. 
\end{theorem}
\begin{theorem}
\label{{Le_1n_Th}}
For any positive integer $n$, the value $\Le_{1, n} (-1)$ belongs to the ring 
$\Q[\ln2, \zeta(2), \dots, \zeta(n+1)]$, namely, the following formula
is valid:
$$
2 \Le_{1, n}(-1) =
\sum_{k=1}^n \zeta(\overline k) \zeta(\overline {n-k+1}) - n \zeta(n+1).
$$
\end{theorem}
\proof
If $n$ is even, then the required formula immediately follows from
Theorem \ref{Le_mn_Th} (substitute $m=1$). If $n$ is odd, then use
(\ref{eq7}) with $z=-1$:
$$
2 \Li_{1,n}(-1)=\sum_{k=1}^n (-1)^{k+1}
\zeta(\overline k) \zeta(\overline {n+1-k}).
$$
Hence
\begin{align*}
2 \Le_{1,n}(-1) & =2 \Li_{1,n}(-1) + 2\zeta(\overline {n+1}) \\
& =
\sum_{k=1}^n \zeta(\overline k) \zeta(\overline {n+1-k})
- 2 \sum_{p=0}^{(n+1)/2} \zeta(\overline {2p}) \zeta(\overline {n+1-2p}).
\end{align*}
We transform the subtrahend into $n \zeta(n+1)$. It is known that
(\cite[1.12(2), 1.13(22)]{Bateman})
\begin{equation*}
\zeta(\overline {2p})=
(-1)^p (1-2^{1-2p}) \frac{(2 \pi)^{2p}}{2 (2p)!} B_{2p}, \quad
\zeta(2q)=
(-1)^{q+1} \frac{(2 \pi)^{2q}}{2 (2q)!} B_{2q},
\end{equation*}
for integer $p \ge 0$ and $q \ge 1$.
Therefore the required equality
$$
2 \sum_{p=0}^{(n+1)/2} \zeta(\overline {2p}) \zeta(\overline {n+1-2p})=
n \zeta(n+1)
$$
is equivalent to the following relation for the Bernoulli numbers:
$$
\sum_{p=0}^q (1-2^{1-2p}) (1-2^{1-(2q-2p)})
\frac{B_{2p} B_{2q-2p}}{(2p)! (2q-2p)!} =
-(2q-1) \frac{B_{2q}}{(2q)!}.
$$
This is a special case of W.~Gosper's equality
(see \cite[(34)]{wolfram})
\begin{equation}
\label{GosperId}
\sum_{p=0}^r (1-2^{1-p}) (1-2^{1-(r-p)})
\frac{B_{p} B_{r-p}}{p! (r-p)!} =
-(r-1) \frac{B_{r}}{r!}.
\end{equation}
Indeed $B_{2k-1}=0$ when $k>1$
(see \cite[1.13(17)]{Bateman}), so $(1-2^{1-k}) B_k = 0$ for
any odd $k \ge 1$.
The proof of (\ref{GosperId}) was kindly communicated to the author
by W.~Gosper. Consider the generating function for the Bernoulli
polynomials (see \cite[1.13(2)]{Bateman}):
$$
\frac{z e^{xz}}{e^z-1} = \sum_{k=0}^{\infty} B_k(x) \frac{z^k}{k!},
\quad |z|<2\pi.
$$
Multiplying the generating functions, we obtain the identity
$$
\sum_{p=0}^{r} \frac{ B_{p}(x) B_{r-p}(y) } {p! (r-p)!} =
\frac{ r(x+y-1) B_{r-1}(x+y-1) - (r-1) B_{r}(x+y-1) }{r!}.
$$
Setting $x=y=1/2$ and using the value $B_k(1/2)=-(1-2^{1-k})B_k$
(see \cite[(23.1.21)]{Abramowitz}), we get
(\ref{GosperId}). The proof of the theorem is finished. \qed

Using
$$
\Li_{m,n}(-1) = \Le_{m,n}(-1)-\Li_{m+n}(-1),
$$
it is easy to write out formulas for the values $\Li_{m,n}(-1)$ 
for odd $m+n$ and $\Li_{1,n}(-1)$ for any $n$.
In particular, the equality
$$
2 \Li_{1, n}(-1) =
\sum_{k=1}^n \zeta(\overline k) \zeta(\overline {n-k+1}) - n \zeta(n+1)
-2 \zeta(\overline {n+1})
$$
is equivalent to Theorem \ref{{Le_1n_Th}} and was proved by N.~Nielsen
(\cite[p. 50, (6)]{Nielsen}).
Theorem \ref{dualityTh} yields
$$
\Le_{\inbr{1}_m, 2, \inbr{1}_n} \left( \frac{1}{2} \right)=
-\Le_{m+1,n+1}(-1),
$$
whence we obtain the values
$\Le_{2, \inbr{1}_n} (1/2)$ for any $n$ and
$\Le_{\inbr{1}_m, 2, \inbr{1}_n} (1/2)$ for odd $m+n$.

In Section \ref{LowWeightSection}, having found all values
of generalized polylogarithms of weight 5, we will obtain 
the curious formula
$$
\Le_{1,3,1}(1/2)=-\Le_{2,1,2}(-1)=
\frac{53}{64} \zeta(5) - \frac{1}{16} \zeta(2) \zeta(3).
$$
An interesting question is whether it is possible to generalize
this formula.

\section{Generalized polylogarithms of small weight}
\label{LowWeightSection}

Suppose the values $\Li_{\vec{s}}(1/2)$ for all vectors $\vec{s}$
of a certain weight have been calculated. Then the associated
values $\Le_{\vec{s}}(1/2)$, $\Le_{\vec{s}}(-1)$, $\Le_{\vec{s}}(-1)$ 
for the same vectors 
can be computed using the equality (\ref{eq5}) and Theorems
\ref{dualityTh} and \ref{LiTransformation}. In this section we study the values
$\Li_{\vec{s}}(1/2)$ for all vectors of weight $\le 5$, 
and in the Appendix we give tables for these values and the associated values.

{\bf Weight 1.} There is only one (generalized) polylogarithm of weight 1, 
and $\Li_1(1/2)=\ln 2$ by equality (\ref{eq1}).

{\bf Weight 2.} The value of $\Li_{1,1}(1/2)$ follows from (\ref{eq6}),
and the value of $\Li_2(1/2)$ is given in (\ref{eq3}).

{\bf Weight 3.} The value of $\Li_{3}(1/2)$ is given in (\ref{eq4}), those of
$\Li_{2,1}(1/2)$ and $\Li_{1,2}(1/2)$ follow from Theorem \ref{Li121_Th}, and
$\Li_{1,1,1}(1/2)$ from (\ref{eq6}).

{\bf Weight 4.} The value 
$\Li_{4}(1/2)$ is considered as a new constant,
$\Li_{1,1,1,1}(1/2)$ is obtained from (\ref{eq6}),
$\Li_{2,1,1}(1/2)$, $\Li_{1,2,1}(1/2)$ and $\Li_{1,1,2}(1/2)$
follow from Theorem \ref{Li121_Th}, $\Li_{2,2}(1/2)$ does from Theorem
\ref{Li2_Th}, $\Li_{1,3}(1/2)$ does from Corollary \ref{cor1}, 
$\Li_{3,1}(1/2)$ does from Corollary \ref{cor3}
(taking into account the found $\Li_{2,2}(1/2)$).
Alternatively, $\Li_{1,3}(1/2)$ and $\Li_{3,1}(1/2)$ can be found 
from formulas given in \cite{BradleyResolutions}.

{\bf Weight 5.} The value
$\Li_{5}(1/2)$ is considered as a new constant;
$\Li_{1,1,1,1,1}(1/2)$ is obtained from (\ref{eq6});
$\Li_{2,1,1,1}(1/2)$, $\Li_{1,2,1,1}(1/2)$, $\Li_{1,1,2,1}(1/2)$ and
$\Li_{1,1,1,2}(1/2)$ follow from Theorem \ref{Li121_Th};
$\Li_{1,2,2}(1/2)$ from Theorem \ref{Li2_Th}; and
$\Li_{1,1,3}(1/2)$ and $\Li_{1,3,1}(1/2)$ from formulas given in
\cite{BradleyResolutions}.
It is left to find the seven values of $\Li$ with indices
$$
(4,1), (3,2), (3,1,1), (2,3), (2,2,1), (2,1,2), (1,4).
$$
Theorem \ref{th2} yields the sum of those seven values
(in terms of values already found); it is equal to the dot product
\begin{equation}
\label{Weight5Eq1}
\left( \frac{59}{32}, \frac{1}{16},
-\frac{7}{8}, 0, 0, \frac{1}{120}, -2, 0 \right) \cdot \vec{B}, 
\end{equation}
where
$$
\vec{B} = 
( \zeta(5), \zeta(2) \zeta(3), \zeta(4) \ln 2, \zeta(3) \ln^2 2,
\zeta(2) \ln^3 2, \ln^5 2, \Li_5(1/2), \Li_4(1/2) \ln 2).
$$
Theorem \ref{HolderTh} with $z=1/2$ and
$y_1 y_2 y_3 y_4 y_5 = x_0^3 x_1^2$ gives
\begin{align*}
\Li_{4,1}\left( \frac{1}{2} \right) + 
\Li_{3,1,1}\left( \frac{1}{2} \right) = & \;
\zeta(4,1)
- \Li_{3,1} \left( \frac{1}{2} \right) \Li_{1}\left( \frac{1}{2} \right)
- \Li_{2,1} \left( \frac{1}{2} \right) \Li_{1,1}\left( \frac{1}{2} \right) \\
& - \Li_{1,1} \left( \frac{1}{2} \right) \Li_{1,1,1}\left( \frac{1}{2} \right)
- \Li_{1} \left( \frac{1}{2} \right) \Li_{2,1,1}\left( \frac{1}{2} \right).
\end{align*}
Using values of generalized polylogarithms of weight $\le 4$ and the value
$\zeta(4,1)=2 \zeta(5) - \zeta(2) \zeta(3)$
(for evaluation of multiple zeta values of small weight
and tables of them, see for example \cite{MinhPetitot}), we obtain
\begin{equation}
\label{Weight5Eq2}
\Li_{4,1}\left( \frac{1}{2} \right) + 
\Li_{3,1,1}\left( \frac{1}{2} \right) = 
\left( 2, -1, -\frac{9}{8}, \frac{15}{16},
-\frac{1}{4}, \frac{1}{24}, 0, 1 \right) \cdot \vec{B}, 
\end{equation}
Similarly using
$$
\zeta(3,2)=-\frac{11}{2} \zeta(5) + 3 \zeta(2) \zeta(3),
\quad
\zeta(2,3)=\frac{9}{2} \zeta(5) - 2 \zeta(2) \zeta(3)
$$
we find
\begin{equation}
\label{Weight5Eq3}
\Li_{3,2}\left( \frac{1}{2} \right) + 
\Li_{2,2,1}\left( \frac{1}{2} \right) = 
\left( -\frac{11}{2}, \frac{47}{16}, \frac{47}{16}, -\frac{45}{16},
\frac{5}{6}, -\frac{1}{8}, 0, -3 \right) \cdot \vec{B}, 
\end{equation}
\begin{equation}
\label{Weight5Eq4}
\Li_{2,3}\left( \frac{1}{2} \right) + 
\Li_{2,1,2}\left( \frac{1}{2} \right) = 
\left( \frac{9}{2}, -\frac{37}{16}, -\frac{43}{16}, \frac{37}{16},
-\frac{3}{4}, \frac{1}{8}, 0, 3 \right) \cdot \vec{B}.
\end{equation}
Summing (\ref{Weight5Eq2}), (\ref{Weight5Eq3}), (\ref{Weight5Eq4}) 
and subtracting the result from (\ref{Weight5Eq1}), we find that
$$
\Li_{1,4}\left( \frac{1}{2} \right) = 
\left( 
\frac {27}{32}, \frac {7}{16}, 0, -\frac {7}{16},
	\frac {1}{6}, -\frac {1}{30}, -2, -1 \right) \cdot \vec{B}.
$$
From (\ref{ShuffleWeight5Eq1}), (\ref{ShuffleWeight5Eq2}),
(\ref{ShuffleWeight5Eq3}) with $z=1/2$ we obtain
\begin{equation}
\label{Weight5Eq5}
2 \Li_{4,1}\left( \frac{1}{2} \right) + 
\Li_{3,2}\left( \frac{1}{2} \right) +
\Li_{2,3}\left( \frac{1}{2} \right) = 
\left( -\frac{27}{32}, -\frac{7}{16}, 0, \frac{7}{16},
	-\frac{1}{6}, \frac{1}{30}, 2, 2 \right) \cdot \vec{B},
\end{equation}
\begin{equation}
\label{Weight5Eq6}
6 \Li_{4,1}\left( \frac{1}{2} \right) + 
3 \Li_{3,2}\left( \frac{1}{2} \right) +
\Li_{2,3}\left( \frac{1}{2} \right) = 
\left( 0, \frac{7}{16}, -\frac{5}{8}, -\frac{7}{16},
	\frac{1}{3}, -\frac{1}{12}, 0, 0 \right) \cdot \vec{B},
\end{equation}
\begin{equation}
\label{Weight5Eq7}
2 \Li_{2,2,1}\left( \frac{1}{2} \right) + 
2 \Li_{2,1,2}\left( \frac{1}{2} \right) = 
\left( -\frac{3}{16}, \frac{1}{8}, 0, \frac{1}{8},
	-\frac{1}{6}, \frac{1}{30}, 0, 0 \right) \cdot \vec{B}.
\end{equation}
From the six linearly independent equations (\ref{Weight5Eq2})-(\ref{Weight5Eq7})
we find six independent values.

All values of generalized polylogarithms of weight
$\le 5$ at the points $z = 1/2$ and $z = -1$ are given in the Appendix.

Unfortunately, the shuffle relations and equalities from Section \ref{HalfSection}
do not allow us to find all values of generalized polylogarithms
of a fixed weight $\ge 6$ at the points $z = 1/2$ and $z = -1$.
(In the next section we discuss
the dimension of the respective linear space and the introduction of 
new constants.)

Some values of generalized polylogarithms of small weight were obtained long ago. For example, the equality
$$
\Le_{2,1}\left( \frac{1}{2} \right)=\zeta(3)-\frac{1}{2} \zeta(2) \ln 2
$$
was found by S.~Ramanujan (see \cite[p. 258]{Berndt}).
Turning now to alternating multiple zeta values, in \cite{BroadhurstEnumeration}, \cite{BroadhurstConjectured}
D.J.~Broadhurst obtained a series of relations for them and advanced
many interesting conjectures. He not only used numerical experimentation
with the help the PSLQ algorithm, but also rigorously found values 
of certain length and weight (for example, double series
of weight $\le 44$). Another approach to studying values of 
small weight was proposed by J.A.M.~Vermaseren
(\cite[Section 5]{Vermaseren}). He calculated all values up to weight 9
(in the refereed paper, only up to weight 7).
Vermaseren's method in fact uses a conjecture explicitly stated 
in \cite[Conjecture 7.1]{Jacob}:
\begin{conjecture}
All algebraic relations among alternating multiple zeta values
are generated by shuffle and stuffle relations.
\end{conjecture}
Assuming this conjecture, M.~Bigotte, an author of \cite{Jacob}, found all 
alternating multiple zeta values of weight $\le 8$.
                  
\def\F{{\mathbb F}}
\def\D{{\mathbb D}}
\def\B{{\mathcal{B}}}
\def\L{{\mathcal{L}}}
\section{The space generated by values of generalized polylogarithms}

Let $\L_w$ be the $\Q$-linear space generated by values $\Li_{\vec{s}}(1/2)$
with vectors $\vec{s}$ of weight $w$. With the help of the
linear relations (\ref{eq5}) and Theorem \ref{dualityTh}
it is easy to show that the three linear spaces generated by
$\{ \Le_{\vec{s}}(1/2) \}_{w(\vec{s})=w}$ or
$\{ \Li_{\vec{s}}(-1) \}_{w(\vec{s})=w}$ or
$\{ \Le_{\vec{s}}(-1) \}_{w(\vec{s})=w}$
all coincide with $\L_w$.
If Conjecture \ref{EulerSumsConj} is true, then $\L_w$ also contains all
alternating multiple zeta values of weight $w$.

Denote by $\F_w$ the dimension of $\L_w$ over $\Q$.
Obviously $\F_1=1$. From Section \ref{LowWeightSection} it follows that
$\F_2 \le 2$, $\F_3 \le 3$ and $\F_4 \le 5$.
Moreover $\F_2=2$ by virtue of the linear independence of
$\pi^2$ and $\ln^2 2$ over $\Q$, which follows from
the transcendence of the number $\pi/\ln 2$. 
Indeed suppose that $\pi/\ln 2$ is algebraic.
Then the number $\pi i/\ln 2$ is algebraic, but not rational
(because it is purely imaginary).
From Gelfond's theorem we find that the number
$2^{\pi i/\ln 2}=e^{\pi i}=-1$ is transcendental, a contradiction.

\begin{conjecture}
\label{conj1}
For $w \ge 1$ we have $\F_w=f_w$, where $\{ f_w \}$ is the Fibonacci sequence,
defined by the recurrence relation $f_w=f_{w-2}+f_{w-1}$ and
initial values $f_0=f_1=1$.
\end{conjecture}

Denote by $\F_w^{A}$ the dimension over $\Q$ of the $\Q$-linear space generated
by values of the alternating multiple zeta function of weight $w$. 
Since this space contains $\{\Li_{\vec{s}}(-1)\}_{w(\vec{s})=w}$,
we have $\F_w^{A} \ge \F_w$.
If Conjecture \ref{EulerSumsConj} is valid, then $\F_w^{A} = \F_w$.
In \cite[(14)]{BroadhurstEnumeration} the equality
$\F_w^{A} = f_w$ was conjectured, while in \cite{DeligneGoncharov}
the estimate $\F_w^{A} \le f_w$ was proved (hence $\F_w \le f_w$).

By Conjecture \ref{conj1} we have $\F_6=13$. However, there are only 12
"classical" constants of weight 6 that are involved in expressions for values of generalized polylogarithms
and which are (presumably) linearly independent over $\Q$, namely,
$$
\zeta(6), \zeta(3)^2, \zeta(5) \ln 2, \zeta(2) \zeta(3) \ln 2, \zeta(4) \ln^2 2,
\zeta(3) \ln^3 2, \zeta(2) \ln^4 2,
$$
$$
\ln^6 2, \Li_6(1/2), \Li_5(1/2) \ln 2, \Li_4(1/2) \ln^2 2, \Li_4(1/2) \zeta(2).
$$
Therefore (conjecturally) there exists a generalized polylogarithm value
of weight 6 that is not a rational linear combination of 
those 12 constants. Seemingly $\Le_{5,1}(-1)$ is such a value;
in \cite[Section 7]{flajolet} it is written that expressions
for values $\Le_{2n+1,1}(-1)$ (denoted by $- \mu_n$ in the paper), $n \ge 2$,
in terms of classical constants are unknown.
Another example is the value $\Li_{2,2,1,1}(1/2)$. 
In both cases we failed to find numerically an expression
in terms of the 12 constants mentioned above, but we succeeded 
for the combination $\Li_{2,2,1,1}(1/2)+(9/4) \Le_{5,1}(-1)$.

In the work \cite{BroadhurstEnumeration}, instead of $\Le_{5,1}(-1)$
the constant $\zeta(\overline 5, \overline 1)$ was proposed for weight 6, 
and for weight 7 two new constants are needed, for example
$\zeta(5, 1, \overline 1)$ and $\zeta(3, 3, \overline 1)$.
In \cite{Vermaseren} the following constants (of weights 6 and 7) are given:
$$
\sum_{n_1 \ge n_2 \ge 1} \frac{(-1)^{n_1+n_2}}{n_1^5 n_2},
\quad
\sum_{n_1 \ge n_2 \ge n_3 \ge 1}
\frac{(-1)^{n_1}}{n_1^5 n_2 n_3},
\quad
\sum_{n_1 \ge n_2 \ge n_3 \ge 1}
\frac{(-1)^{n_2+n_3}}{n_1^5 n_2 n_3}.
$$

The numbers $f_w$ can also be defined by the generating function
$$
\sum_{w=0}^{\infty} f_w x^w = \frac{1}{1-x-x^2}.
$$
It is readily shown (with the help of the generating function or
directly) that $f_w = {\rm Card}(\B_w)$, where
$$
\B_w = 
\left\{
(s_1, s_2, \dots, s_l) :
\sum_{j=1}^l s_j = w, \; s_j \in \{ 1, 2 \}
\right\}.
$$
It is natural to advance the following conjecture.
\begin{conjecture}
As a basis of $\L_w$ it is possible to take
any of the following sets of numbers: \\
1) $\{ \Li_{\vec{s}}(1/2) : \vec{s} \in \B_w \}$ \\
2) $\{ \Li_{\vec{s}}(-1) : \vec{s} \in \B_w \}$ \\
3) $\{ \Le_{\vec{s}}(-1) : \vec{s} \in \B_w \}$.
\end{conjecture}
We don't include here the variant with the function $\Le_{\vec{s}}(z)$
at the point $z=1/2$ because it would fail due to the equality
$$
5 \Le_{1,1,1} \left( \frac{1}{2} \right)=
6 \Le_{1,2} \left( \frac{1}{2} \right).
$$
The following proposition shows that the two variants 1) and 2)
with $\Li_{\vec{s}}(z)$ are in fact equivalent.
\begin{claim}
For any positive integer $w$, the two linear spaces over $\Q$ generated
by the sets of numbers $\Li_{\vec{s}}(1/2)$, ${\vec{s} \in \B_w}$, or
$\Li_{\vec{t}}(-1)$, ${\vec{t} \in \B_w}$, coincide.
In particular, the linear independence over $\Q$ of the numbers
in one set implies the linear independence of the numbers in the other set.
\end{claim}
\proof
Theorem \ref{LiTransformation} yields that for $\vec{s} \in \B_w$ 
the function $\Li_{\vec{s}}(\frac{-z}{1-z})$ is represented as an
integer linear combination of $\Li_{\vec{t}}(z)$, $\vec{t} \in B_w$.
Applying that fact for $z=-1$ and $z=1/2$, we obtain that
each number $\Li_{\vec{s}}(1/2)$, ${\vec{s} \in \B_w}$, is an 
integer linear combination of numbers
$\Li_{\vec{t}}(-1)$, ${\vec{t} \in \B_w}$, and vice versa. 
That proves the proposition. \qed

We now consider the variant 1) with the function $\Li_{\vec{s}}(z)$ 
at the point $z=1/2$. For that we divide the basis conjecture
into two independent conjectures.

\begin{conjecture}
\label{conj2}
The numbers $\Li_{\vec{s}}(1/2)$, $\vec{s} \in \B_w$, are linearly 
independent over $\Q$.
\end{conjecture}
\begin{conjecture}
\label{conj3}
For any vector $\vec{s_0}$ of weight $w$ the value $\Li_{\vec{s_0}}(1/2)$
is a rational linear combination of the numbers 
$\Li_{\vec{s}}(1/2)$, $\vec{s} \in \B_w$.
\end{conjecture}
From Conjecture \ref{conj2} it follows that $\F_w \ge f_w$, 
while Conjecture \ref{conj3} gives the inequality $\F_w \le f_w$.
Conjecture \ref{conj2} is true at least for $w=1$ and $w=2$. 
Conjecture \ref{conj3} for $w \le 5$ is readily verified
with the help of tables given in the Appendix.
For $w \le 9$ this conjecture was checked numerically
with the help of PSLQ, that is, the required expressions were found and
hold with high precision.

Note the similarity with the dimension $\D_w$ of the space 
generated by multiple zeta values $\zeta(s_1, \dots, s_l)$ of weight $w$.
D.~Zagier (\cite{zagier}) conjectured the following. 
\begin{conjecture}
For $w \ge 2$ we have $\D_w=d_w$, where the sequence $\{ d_w \}$
is defined by the recurrence relation
$d_w=d_{w-3}+d_{w-2}$ with initial values $d_0=1$, $d_1=0$, $d_2=1$.
\end{conjecture}
The numbers $d_w$ can be defined with the help of the generating function
$$
\sum_{w=0}^{\infty} d_w x^w = \frac{1}{1-x^2-x^3}.
$$
In \cite{DeligneGoncharov}, \cite{terasoma} the inequality 
$\D_w \le d_w$ is proved. No nontrivial lower estimate for $\D_w$
has been proved (here the estimate $\D_w \ge 1$ is called trivial). 

M.E.~Hoffman (\cite{hoffman}) advanced the following consistent conjecture.
\begin{conjecture}
\label{HoffmanConj}
As a basis of the $\Q$-linear space generated by the values $\zeta(\vec{s})$
of weight $w$ it is possible to take the numbers
$\zeta(t_1, t_2, \dots, t_l)$ with
$\sum_{j=1}^l t_j = w$ and $t_j \in \{ 2, 3 \}$.
\end{conjecture}
Hoang Ngoc Minh has checked that for any vector
$\vec{s_0}$ of weight $\le 16$ the value
$\zeta(\vec{s_0})$ is a rational linear combination of multiple zeta values
in the presumptive basis.

If $\vec{s}=(s_1, s_2, \dots, s_l)$, then
$\Li_{\vec{s}}(-1)=\zeta(\overline s_1, s_2, \dots, s_l)$. Hence
if Conjectures \ref{EulerSumsConj}, \ref{conj2}, and \ref{conj3} are true,
then the following conjecture is also true.
\begin{conjecture}
As a basis of the $\Q$-linear space generated by the values of
the alternating multiple zeta function $\zeta(\vec{s}; \vec{\sigma})$
of weight $w$ it is possible to take the numbers
$\zeta(\overline t_1, t_2, \dots, t_l)$ with
$\sum_{j=1}^l t_j = w$ and $t_j \in \{ 1, 2 \}$.
\end{conjecture}

In conclusion, we note the result of V.N.~Sorokin 
\cite{sorokinPols} on the linear independence 
over $\Q$ of values of generalized polylogarithms
at a rational point $p/q \ne 0$ lying near zero.
Namely, the generalized polylogarithms $\Li_{\vec{s}}(z)$
of weight $w(\vec{s})\le r$ are linearly independent
at a point $z=p/q$ if $|p/q| < (2pe^r)^{-N_r}$, where $N_r=2^r-1$.

\newpage

\section{Acknowledgments}
The author thanks D.M.~Bradley for consultations on the subject of this paper,
P.~Lison\v ek for providing source code to compute values of 
generalized polylogarithms, and J.~Sondow for help in translating the paper
from Russian and for several suggestions.


\newcommand{\namefont}{\scshape}
\newcommand{\titlefont}{\itshape}
\newcommand{\nomer}{{No. }}


\section*{Appendix}
\addcontentsline{toc}{section}{Appendix}

The following tables give all values of generalized polylogarithms of weight
$\le 5$ at the points $z = 1/2$ and $z = -1$.
The values were found in Section \ref{LowWeightSection}
and do not rely on any conjectures.

\newpage

Values for weight 2 are of the form $A \zeta(2) + B \ln^2 2$
(the table gives the vector $(A,B)$).
\begin{center}
\begin{tabular}
{| c| c| c | c | c |}
\hline
$\vec{s}$ & $\Li_{\vec{s}}(1/2)$ & $\Le_{\vec{s}}(1/2)$
	& $\Li_{\vec{s}}(-1)$ & $\Le_{\vec{s}}(-1)$ \\
\hline
2 & $(\frac{1}{2},-\frac{1}{2})$ & $(\frac{1}{2},-\frac{1}{2})$
	& $(-\frac{1}{2},0)$ & $(-\frac{1}{2},0)$ \\
\hline
1,1 & $(0,\frac{1}{2})$ & $(\frac{1}{2},0)$
	& $(0,\frac{1}{2})$ & $(-\frac{1}{2},\frac{1}{2})$ \\
\hline
\end{tabular}
\end{center}

Values for weight 3 are of the form
$A \zeta(3) + B \zeta(2) \ln 2 + C \ln^3 2$
(the table gives the vector $(A,B,C)$).
\begin{center}
\begin{tabular}
{| c| c| c | c | c |}
\hline
$\vec{s}$ & $\Li_{\vec{s}}(1/2)$ & $\Le_{\vec{s}}(1/2)$
	& $\Li_{\vec{s}}(-1)$ & $\Le_{\vec{s}}(-1)$ \\
\hline
3 & $(\frac{7}{8},-\frac{1}{2},\frac{1}{6})$
	& $(\frac{7}{8},-\frac{1}{2},\frac{1}{6})$
	& $(-\frac{3}{4},0,0)$
	& $(-\frac{3}{4},0,0)$ \\
\hline
2,1 & $(\frac{1}{8},0,-\frac{1}{6})$
	& $(1,-\frac{1}{2},0)$
	& $(\frac{1}{8},0,0)$
	& $(-\frac{5}{8},0,0)$ \\
\hline
1,2 & $(-\frac{1}{4},\frac{1}{2},-\frac{1}{6})$
	& $(\frac{5}{8},0,0)$
	& $(-\frac{1}{4},\frac{1}{2},0)$
	& $(-1,\frac{1}{2},0)$ \\
\hline
1,1,1 & $(0,0,\frac{1}{6})$
	& $(\frac{3}{4},0,0)$
	& $(0,0,-\frac{1}{6})$
	& $(-\frac{7}{8},\frac{1}{2},-\frac{1}{6})$ \\
\hline
\end{tabular}
\end{center}

Values for weight 4 are of the form $A \zeta(4) + B \zeta(3) \ln 2
+ C \zeta(2) \ln^2 2 + D \ln^4 2 + E \Li_4(1/2)$
(the table gives the vector  $(A,B,C,D,E)$).

\begin{center}
\noindent
{\small
\begin{tabular}
{| c| c| c | c | c |}
\hline
$\vec{s}$ & $\Li_{\vec{s}}(1/2)$ & $\Le_{\vec{s}}(1/2)$
	& $\Li_{\vec{s}}(-1)$ & $\Le_{\vec{s}}(-1)$ \\
\hline
4 & $(0,0,0,0,1)$
	& $(0,0,0,0,1)$
	& $(-\frac{7}{8},0,0,0,0)$
	& $(-\frac{7}{8},0,0,0,0)$ \\
\hline
3,1 & $(\frac{1}{8},-\frac{1}{8},0,\frac{1}{24},0)$
	& $(\frac{1}{8},-\frac{1}{8},0,\frac{1}{24},1)$
	& $(-\frac{15}{8},\frac{7}{4},-\frac{1}{2},\frac{1}{12},2)$
	& $(-\frac{11}{4},\frac{7}{4},-\frac{1}{2},\frac{1}{12},2)$ \\
\hline
2,2 & $(\frac{1}{16},\frac{1}{4},-\frac{1}{4},\frac{1}{24},0)$
	& $(\frac{1}{16},\frac{1}{4},-\frac{1}{4},\frac{1}{24},1)$
	& $(\frac{65}{16},-\frac{7}{2},1,-\frac{1}{6},-4)$
	& $(\frac{51}{16},-\frac{7}{2},1,-\frac{1}{6},-4)$ \\
\hline
1,3 & $(-\frac{5}{16},\frac{7}{8},-\frac{1}{4},\frac{1}{24},0)$
	& $(-\frac{5}{16},\frac{7}{8},-\frac{1}{4},\frac{1}{24},1)$
	& $(-\frac{5}{16},\frac{3}{4},0,0,0)$
	& $(-\frac{19}{16},\frac{3}{4},0,0,0)$ \\
\hline
2,1,1 & $(1,-\frac{7}{8},\frac{1}{4},-\frac{1}{12},-1)$
	& $(\frac{19}{16},-\frac{3}{4},0,0,0)$
	& $(-1,\frac{7}{8},-\frac{1}{4},\frac{1}{24},1)$
	& $(\frac{5}{16},-\frac{7}{8},\frac{1}{4},-\frac{1}{24},-1)$ \\
\hline
1,2,1 & $(-3,\frac{11}{4},-\frac{3}{4},\frac{1}{12},3)$
	& $(-\frac{51}{16},\frac{7}{2},-1,\frac{1}{6},4)$
	& $(3,-\frac{11}{4},\frac{3}{4},-\frac{1}{8},-3)$
	& $(-\frac{1}{16},-\frac{1}{4},\frac{1}{4},-\frac{1}{24},-1)$ \\
\hline
1,1,2 & $(3,-\frac{23}{8},1,-\frac{1}{6},-3)$
	& $(\frac{11}{4},-\frac{7}{4},\frac{1}{2},-\frac{1}{12},-2)$
	& $(-3,\frac{23}{8},-1,\frac{1}{8},3)$
	& $(-\frac{1}{8},\frac{1}{8},0,-\frac{1}{24},-1)$ \\
\hline
1,1,1,1 & $(0,0,0,\frac{1}{24},0)$
	& $(\frac{7}{8},0,0,0,0)$
	& $(0,0,0,\frac{1}{24},0)$
	& $(0,0,0,0,-1)$ \\
\hline
\end{tabular}
}
\end{center}

Values for weight 5 are of the form 
$$
A \zeta(5) + B \zeta(2) \zeta(3)
+ C \zeta(4) \ln 2 + D \zeta(3) \ln^2 2 
$$
$$
+\ E \zeta(2) \ln^3 2 + F \ln^5 2 + 
G \Li_5(1/2) + H \Li_4(1/2) \ln 2
$$
(the tables give the vector $(A,B,C,D,E,F,G,H)$).
\begin{center}
{\small
\begin{tabular}
{| c| c| c |}
\hline
$\vec{s}$ & $\Li_{\vec{s}}(1/2)$ & $\Le_{\vec{s}}(1/2)$ \\
\hline
5 & $ (0, 0, 0, 0, 0, 0, 1, 0) $
	& $(0, 0, 0, 0, 0, 0, 1, 0)$ \\
\hline
4,1 & $ (\frac{1}{32}, -\frac{1}{2}, -\frac{1}{8}, \frac {1}{2}, 
	 -\frac{1}{6}, \frac{1}{40}, 1, 1)$
	& $(\frac {1}{32}, -\frac {1}{2}, -\frac {1}{8}, \frac {1}{2},
	-\frac {1}{6}, \frac {1}{40}, 2, 1)$ \\
\hline
3,2 & $ (\frac {23}{64}, \frac {23}{16}, -\frac {1}{16} , -\frac {23}{16},
	\frac {7}{12}, -\frac {13}{120}, -3, -3) $
	& $(\frac {23}{64}, \frac {23}{16}, -\frac {1}{16}, -\frac {23}{16},
	\frac {7}{12}, -\frac {13}{120}, -2, -3)$ \\
\hline
3,1,1 & $ (\frac {63}{32}, -\frac {1}{2}, -1, \frac {7}{16},
	-\frac {1}{12}, \frac {1}{60}, -1, 0) $
	& $(\frac {151}{64}, \frac {7}{16}, -\frac {19}{16}, -\frac {1}{2},
	\frac {1}{3}, -\frac {1}{15}, -2, -2)$ \\
\hline
2,3 & $ (-\frac {81}{64}, -\frac {7}{8}, \frac {5}{16}, \frac {7}{8},
	-\frac {5}{12}, \frac {11}{120}, 3, 3) $
	& $(-\frac {81}{64}, -\frac {7}{8}, \frac {5}{16}, \frac {7}{8},
	-\frac {5}{12}, \frac {11}{120}, 4, 3)$ \\
\hline
2,2,1 & $ (-\frac {375}{64}, \frac {3}{2}, 3, -\frac {11}{8},
	\frac {1}{4}, -\frac {1}{60}, 3, 0) $
	& $(-\frac {227}{32}, \frac {1}{8}, \frac {51}{16}, 0,
	-\frac {1}{3}, \frac {1}{10}, 8, 4)$ \\
\hline
2,1,2 & $ (\frac {369}{64}, -\frac {23}{16}, -3, \frac {23}{16},
	-\frac {1}{3}, \frac {1}{30}, -3, 0) $
	& $(\frac {311}{64}, -\frac {7}{8}, -\frac {11}{4}, \frac {7}{8},
	-\frac {1}{6}, \frac {1}{60}, -2, 0)$ \\
\hline
2,1,1,1 & $ (1, 0, 0, -\frac {7}{16},
	\frac {1}{6}, -\frac {1}{24}, -1, -1) $
	& $(2, -\frac {3}{8}, -\frac {7}{8}, 0,
	0, 0, 0, 0)$ \\
\hline
1,4 & $ (\frac {27}{32}, \frac {7}{16}, 0, -\frac {7}{16},
	\frac {1}{6}, -\frac {1}{30}, -2, -1) $
	& $(\frac {27}{32}, \frac {7}{16}, 0, -\frac {7}{16},
	\frac {1}{6}, -\frac {1}{30}, -1, -1)$ \\
\hline
1,3,1 & $ (-\frac {3}{64}, 0, \frac {1}{8}, -\frac {1}{16},
	0, \frac {1}{120}, 0, 0) $
	& $(\frac {53}{64}, -\frac {1}{16}, 0, 0,
	0, 0, 0, 0)$ \\
\hline
1,2,2 & $ (\frac {3}{16}, -\frac {1}{8}, \frac {1}{16}, \frac {1}{8},
	-\frac {1}{12}, \frac {1}{120}, 0, 0) $
	& $(\frac {89}{64}, \frac {7}{4}, 0, -\frac {7}{4},
	\frac {2}{3}, -\frac {2}{15}, -4, -4)$ \\
\hline
1,2,1,1 & $ (-4, 0, 1, \frac {7}{8}, 
	-\frac {5}{12}, \frac {1}{12} , 4, 3) $
	& $(-\frac {21}{32}, \frac {3}{4}, 0, 0,
	0, 0, 0, 0)$ \\
\hline
1,1,3 & $ (-\frac {3}{64}, \frac {1}{16}, -\frac {5}{16}, \frac {7}{16},
	-\frac {1}{12}, \frac {1}{120}, 0, 0) $
	& $(-\frac {15}{32},  -\frac {3}{8}, 0, \frac {7}{8},
	-\frac {1}{3}, \frac {1}{15}, 2, 2)$ \\
\hline
1,1,2,1 & $ (6, 0, -3, \frac {1}{16}, 
	\frac {1}{4}, -\frac {1}{12} , -6, -3) $
	& $(-\frac {11}{32},  \frac {5}{8}, 0, 0, 
	0, 0, 0, 0)$ \\
\hline
1,1,1,2 & $ (-4, 0, 3, -1,
	\frac {1}{6}, 0, 4, 1) $
	& $(\frac {59}{32}, -\frac {1}{2}, 0, 0, 
	0, 0, 0, 0)$ \\
\hline
1,1,1,1,1 & $ (0, 0, 0, 0, 0, \frac {1}{120}, 0, 0) $
	& $(\frac {15}{16}, 0, 0, 0, 
	0, 0, 0, 0)$ \\
\hline
\end{tabular}
}
\vskip 3mm
{\small
\begin{tabular}
{| c| c| c |}
\hline
$\vec{s}$ & $\Li_{\vec{s}}(-1)$ & $\Le_{\vec{s}}(-1)$ \\
\hline
5 & $ (-\frac {15}{16}, 0, 0, 0, 0, 0, 0, 0) $
	& $(-\frac {15}{16}, 0, 0, 0, 0, 0, 0, 0)$ \\
\hline
4,1 & $ (-\frac{29}{32}, \frac{1}{2}, 0, 0, 0, 0, 0, 0)$
	& $(-\frac {59}{32}, \frac {1}{2}, 0, 0, 
	0, 0, 0, 0)$ \\
\hline
3,2 & $ (\frac{41}{32}, -\frac{5}{8}, 0, 0, 0, 0, 0, 0) $
	& $(\frac {11}{32},  -\frac {5}{8}, 0, 0, 
	0, 0, 0, 0)$ \\
\hline
3,1,1 & $ (\frac{33}{32}, \frac{1}{2}, 0, -\frac{7}{8},
	\frac{1}{3}, -\frac{1}{15}, -2, -2) $
	& $(\frac {15}{32},  \frac {3}{8}, 0, -\frac {7}{8},
	\frac {1}{3}, -\frac {1}{15}, -2, -2)$ \\
\hline
2,3 & $ (\frac{51}{32}, -\frac{3}{4}, 0, 0, 0, 0, 0, 0) $
	& $(\frac {21}{32}, -\frac {3}{4}, 0, 0,
	0, 0, 0, 0)$ \\
\hline
2,2,1 & $ (-\frac{73}{64}, -\frac{3}{2}, 0, \frac{7}{4}, -\frac{2}{3},
	\frac{2}{15}, 4, 4) $
	& $(-\frac {89}{64}, -\frac {7}{4}, 0, \frac {7}{4},
	-\frac {2}{3}, \frac {2}{15}, 4, 4)$ \\
\hline
2,1,2 & $ (-\frac{177}{64}, \frac{23}{16}, 0, 0, 0, 0, 0, 0) $
	& $(-\frac {53}{64}, \frac {1}{16}, 0, 0,
	0, 0, 0, 0)$ \\
\hline
2,1,1,1 & $ (1, 0, 0, -\frac{7}{16}, \frac{1}{6}, -\frac{1}{30}, -1, -1) $
	& $(-\frac {27}{32}, -\frac {7}{16}, 0, \frac {7}{16},
	-\frac {1}{6}, \frac {1}{30}, 1, 1)$ \\
\hline
1,4 & $ (-\frac{17}{16}, \frac{3}{8}, \frac{7}{8}, 0, 0, 0, 0, 0) $
	& $(-2, \frac {3}{8}, \frac {7}{8}, 0,
	0, 0, 0, 0)$ \\
\hline
1,3,1 & $ (-\frac{125}{64}, 0, \frac{15}{8}, -\frac{7}{8},
	\frac{1}{6}, -\frac{1}{60}, 2, 0) $
	& $(-\frac {311}{64}, \frac {7}{8}, \frac {11}{4}, -\frac {7}{8},
	\frac {1}{6}, -\frac {1}{60}, 2, 0)$ \\
\hline
1,2,2 & $ (\frac{125}{16}, \frac{1}{8}, -\frac{65}{16}, 0,
	\frac{1}{3}, -\frac{1}{10}, -8, -4) $
	& $(\frac {227}{32}, -\frac {1}{8}, -\frac {51}{16}, 0,
	\frac {1}{3}, -\frac {1}{10}, -8, -4)$ \\
\hline
1,2,1,1 & $ (-4, 0, 1, \frac{7}{8}, -\frac{5}{12}, \frac{11}{120}, 4, 3) $
	& $(\frac {81}{64}, \frac {7}{8}, -\frac {5}{16}, -\frac {7}{8},
	\frac {5}{12}, -\frac {11}{120}, -4, -3)$ \\
\hline
1,1,3 & $ (-\frac{125}{64}, -\frac{1}{16}, \frac{5}{16}, \frac{1}{2},
	-\frac{1}{3}, \frac{1}{15}, 2, 2) $
	& $(-\frac {151}{64}, -\frac {7}{16}, \frac {19}{16}, \frac {1}{2},
	-\frac {1}{3}, \frac {1}{15}, 2, 2)$ \\
\hline
1,1,2,1 & $ (6, 0, -3, \frac{1}{16}, \frac{1}{4}, -\frac{3}{40}, -6, -3) $
	& $(-\frac {23}{64}, -\frac {23}{16}, \frac {1}{16}, \frac {23}{16},
	-\frac {7}{12}, \frac {13}{120}, 2, 3)$ \\
\hline
1,1,1,2 & $ (-4, 0, 3, -1, \frac{1}{6}, \frac{1}{120}, 4, 1) $
	& $(-\frac {1}{32}, \frac {1}{2}, \frac {1}{8}, -\frac {1}{2},
	\frac {1}{6}, -\frac {1}{40}, -2, -1)$ \\
\hline
1,1,1,1,1 & $ (0, 0, 0, 0, 0, -\frac{1}{120}, 0, 0) $
	& $ (0, 0, 0, 0, 0, 0, -1, 0) $ \\
\hline
\end{tabular}
}
\end{center}

\end {document}